\documentclass[11pt]{article}

\usepackage{amssymb, amsmath, amsthm, graphicx}
\usepackage[left=1in,top=1in,right=1in]{geometry}

\date{}

\theoremstyle{plain}
      \newtheorem{theorem}{Theorem}[section]

      \newtheorem{problem}[theorem]{Problem}
      
      \newtheorem{corollary}[theorem]{Corollary}

\theoremstyle{definition}

\theoremstyle{remark}

\title{Improved bounds for the Ramsey number of tight cycles versus cliques}

\author{Dhruv Mubayi\thanks{Department of Mathematics, Statistics, and Computer Science, University of Illinois, Chicago, IL, 60607 USA.  Research partially supported by NSF grant DMS-1300138. Email: {\tt mubayi@uic.edu}} }

\begin{document}

\maketitle

\begin{abstract}
The 3-uniform tight cycle $C_s^3$ has vertex set $ Z_s$
and edge set $\{\{i, i+1, i+2\}: i \in Z_s\}$.
We prove that for
every $s \not\equiv 0$ (mod 3) and $s \ge 16$ or $s \in \{8,11,14\}$ there is a $c_s>0$ such that the 3-uniform hypergraph Ramsey number $r(C_s^3, K_n^3)$ satisfies
$$r(C_s^3, K_n^3)< 2^{c_s n \log n}.$$ 
This answers in strong form a question of the author and R\"odl who asked for an upper bound of the form
$2^{n^{1+\epsilon_s}}$ for each fixed $s \ge 4$, where $\epsilon_s \rightarrow 0$ as $s \rightarrow \infty$ and $n$ is sufficiently large.  The result is nearly tight as the lower bound is known to be exponential in $n$.

\end{abstract}

\section{Introduction}
A triple system or $3$-graph $H$ with vertex set $V(H)$ is a collection of $3$-element subsets of $V(H)$. Write $K_n^3$ for the complete $3$-graph with vertex set of size $n$ . 
 Given $3$-graphs $F, G$, the Ramsey number $r(F,G)$ is the minimum $n$ such that every red/blue coloring of $K_n^3$ results in a monochromatic red copy of $F$ or a monochromatic blue copy of $G$. In this paper we consider the 3-graph Ramsey number of cycles versus complete 3-graphs.

 For fixed $s \ge 3$ the graph  Ramsey number $r(C_s, K_t)$ has been extensively studied. The case $s=3$ is one of the oldest questions in Ramsey theory and it is known that
$r(C_3,K_t)=\Theta(t^2/\log t)$ (see~\cite{AKS, K} and \cite{BK2, FGM} for recent improvements).
 The next case $r(C_4, K_t)$ seems substantially more difficult and the best known upper and lower bounds are $O(t^2/\log^2t)$ and $\Omega(t^{3/2}/\log t)$, respectively.
An old open problem of Erd\H os asks whether there is a positive $\epsilon$ for which 
$r(C_4, K_t) = O(t^{2-\epsilon})$.  For larger cycles, the best known bounds can be found in~\cite{BK, SU}, and the order of magnitude of $r(C_s, K_t)$ is not known for any fixed $s \ge 4$.

There are several natural ways to define a cycle in hypergraphs. Two possibilities are to consider loose cycles and tight cycles. Here we consider tight cycles.
 For $s > 3$, the tight cycle $C_s^3$ is the $3$-graph
with vertex set $Z_s$ (integers modulo $s$) and edge set 
$$\{\{i, i+1, \ldots, i+k-1\}: i \in Z_s\}.$$  We can view the vertex set of $C_s^3$ as $s$ points on a circle and the edge set as the $s$ circular subintervals each containing $3$ consecutive vertices.

The author and R\"odl~\cite{MR} investigated the hypergraph Ramsey number $r(C_s^3, K_n^3)$ for fixed $s\ge 5$ and large $n$. When $s \equiv 0$ (mod 3) the tight cycle $C_s^3$ is 3-partite, and in this case it is trivial to observe that $r(C_s^3, K_n^3)$ grows like a polynomial in $n$. Determining the growth rate of this polynomial appears to be a very difficult problem and the order of magnitude is not known for any $s>3$. When $s \not\equiv 0$ (mod 3) the Ramsey number is exponential in $n$ as shown in~\cite{MR} below. 

\begin{theorem}{\bf(\cite{MR})} \label{MR}
Fix $s \ge 5$ and $s \not\equiv 0$ (mod 3). There are positive constants $c_1$ and $c_2$ such that
$$2^{c_1n}<r(C^3_s, K^3_n)<2^{c_2n^2\log n}.$$ 
\end{theorem}

Note that when $s=4$, the cycle $C_4^3$ is $K_4^3$ and in this case the lower bound was proved much earlier by Erd\H os and Hajnal~\cite{EH72}, and in fact has been improved to $2^{c_1 n \log n}$ more recently by Conlon-Fox-Sudakov~\cite{CFS}.

As $s$ gets large, the tight cycle $C_s^3$ becomes sparser, so one might expect that $r(C_s^3, K_n^3)$ decreases as a function of $n$ (for fixed $s$). This was asked by R\"odl and the author in~\cite{MR}.

\begin{problem} {\bf (\cite{MR})} \label{prob}
Is the following true?
For each fixed $s \ge 4$ there exists  $\epsilon_s$ such that $\epsilon_s \rightarrow 0$ as $s \rightarrow \infty$ and
$$r(C_s^3, K_n^3) < 2^{n^{1+\epsilon_s}}$$
for all sufficiently large $n$.
\end{problem}

In this short note we give an affirmative answer to this problem by proving the following stronger upper bound.

\begin{theorem} \label{main}
Fix a positive integer $s \not\equiv 0$ (mod 3) such that $s\ge 16$ or $s \in \{8,11,14\}$. There is a positive constant $c$ such that 
$$r(C_s^3, K_n^3) < 2^{c n \log n}.$$
\end{theorem}

Proving a similar result for small values of $s$ remains an open problem, in particular, the cases 
$s\in\{ 4,5,7,10,13\}$.

\section{Proof}
The proof of Theorem~\ref{main} has three components. The first is an observation from~\cite{MR} that  the supersaturation phenomenon from extremal hypergraph theory can be applied to hypergraph Ramsey problems. The crucial new observation here is that we can apply this to families of hypergraphs rather than just individual ones.
The second is a strengthening of the original  proof due to Erd\H os and Hajnal~\cite{EH72} that $r(K_4^3\setminus e, K_n^3) < 2^{cn\log n}$ for some absolute $c>0$, where $K_4^3\setminus e$ is the triple system with four vertices and three edges
(another proof of this upper bound has been given by Conlon-Fox-Sudakov~\cite{CFS}, though that proof is not suitable for our purposes). Finally, the third component is an explicit homomorphism of $C_8^3$ into a particular 3-graph $H_6$ on six vertices.

In order to apply the program above to solve Problem~\ref{prob} (and prove Theorem~\ref{main}), we need to find an $F$ for which $r(F, K_n^3) = 2^{n^{1+o(1)}}$ on the one hand, and $C_s^3$ can be embedded in a blowup of $F$ on the other hand. 
We remark that there are no (non 3-partite) 3-graphs $F$ for which  the order of magnitude of  $\log r(F, K_n^3)$ is known, and the only  (non 3-partite) $F$ for which $r(F, K_n^3) = 2^{n^{1+o(1)}}$ is known is $F=K_4^3\setminus e$. But tight cycles to not embed into a blowup of $K_4^3\setminus e$. We overcome this problem by finding  a (new) family $\cal F$ of 3-graphs for which we can prove $r({\cal F}, K_n^3) = 2^{n^{1+o(1)}}$ and yet embed tight cycles in blowups of each member of $\cal F$.

\subsection{Supersaturation}
Given a hypergraph $H$ and vertex $v \in V(H)$, we say that $w \in V(H)$ is a {\it clone} of $v$ if no edge contains both $v$ and $w$ and for every $e \in H$,
 $v \in e$ if and only if $(e\cup\{w\})\setminus v \in H$.
Given a triple system $F$ and a vertex $v$ in $F$, let $F(v)$ be the triple system obtained from $F$ by replacing $v$ with two clones $v_1,v_2$.  We will use the following result of the author and R\"odl.
\begin{theorem}  {\bf (\cite{MR})} \label{super}
Let $F$ be a triple system with $f$ vertices and $v \in V(F)$. Then  
$$r(F(v), K_n^3)< (r(F, K_n^3))^{f}+f.$$
\end{theorem}

A blowup of $F$ is a a hypergraph obtained by successively applying the cloning operation. 
In particular, if each vertex is cloned $p$ times, then denote the obtained blowup as $F(p)$. By applying Theorem \ref{super} repeatedly, we obtain the following easy corollary. 

\begin{corollary} {\bf (\cite{MR})}  \label{supercor}
Fix a $k$-graph $F$ and an integer $p \ge 2$.
There exists  $c=c(F,p)$ such that 
$$r(F(p), K_n^3) < (r(F, K_n^3))^c.$$
\end{corollary}
 Given a finite family $\cal F$ of 3-graphs, define $r({\cal F}, K_n^3)$ to be the minimum $N$ such that every red/blue coloring of $K_N^3$ results in a red copy of some member of ${\cal F}$ or a blue copy of $K_n^3$.
 Also, for a positive  integer $p$, define ${\cal F}(p)=\{F(p):F \in {\cal F}\}$.  By using the pigeonhole principle, a trivial modification of the proof of Corollary~\ref{supercor} yields the following extension to finite families.
 
\begin{corollary}  \label{supercorf}
Fix a finite family ${\cal F}$ of $k$-graphs  and an integer $p \ge 2$.
There exists  $c=c({\cal F},p)$ such that 
$$r({\cal F}(p), K_n^3) < (r({\cal F}, K_n^3))^c.$$
\end{corollary}
 
\subsection{The triple systems $H_5$ and $H_6$}

Define $H_5$ to be the 3-graph with vertex set 
$\{a,b,u,v,w\}$ and edge set $\{abu, abv, abw, auv, buw\}$ and $H_6$ to be the 3-graph with vertex set 
$\{a,b,u,v,w,x\}$ and edge set $\{abu, abv, abw, abx, auv, bwx\}$.  We will need the following result result of Spencer~\cite{S}. Suppose $H$ is an 
$N$ vertex triple system with average degree at most $D$. Then $H$ has an independent set of size at least $(2/3)N/D^{1/2}$.
%\begin{equation} \label{spen}
%\frac{2}{3} \frac{N}{D^{1/2}}.
%\end{equation}

\begin{theorem} \label{K4H5H6}
There is an absolute positive constant $c$ such that 
$$r(\{K_4^3, H_5, H_6\}, K_n^3) < 2^{cn \log n}.$$
\end{theorem}
\proof
Suppose that $H$ is an $N$ vertex triple system with no independent set of
size $n$.  We will show that $H$ contains a copy of 
some member in $\{K_4^3, H_5, H_6\}$ as long as $N\ge 4^{n}(n!)^2$.  By interpreting the edges of $H$ as the red edges in a red/blue coloring of $K_N^3$, this proves the theorem.

We proceed  by induction on $n$. The result is clearly true for $n=1$. Suppose it holds for $n-1$ and let us show it for $n$. If every pair of vertices in $H$ lie in at most $d=4^{n-1}(n-1)!^2$ edges, then the number of edges of $H$ is at most ${N \choose 2}d/3$ and the average degree of $H$ is at most 
 $$D=\frac{{N \choose 2} d}{N}= \frac{d(N-1)}{2}.$$ Therefore by \cite{S},
$H$ has an independent set of size at least 
$$\frac{2}{3} \frac{N}{D^{1/2}} = \frac{2}{3} \frac{N}{(d(N-1)/2)^{1/2}}
>\frac{2}{3}\frac{N^{1/2}}{d^{1/2}} \ge \frac{2}{3}\frac{2^n n!}{2^{n-1}(n-1)!} > n,$$
a contradiction.
 
We may therefore assume that there is  pair of vertices $a,b$ in $H$ that lie in at least $d+1$ edges.  Consider the set $N(a,b)=\{x: abx\in H\}$. Since $|N(a,b)|>d$, by induction, $N(a,b)$ contains a copy of some member of    $\{K_4^3, H_5, H_6\}$ or an independent set $S$ of size at least $n-1$.
We may assume the latter, otherwise we are done. Now we may use $a$ or $b$ to enlarge the independent set $S$ by one. If we  succeed, then we obtain an independent set of size $n$, which is a contradiction, so we may assume that there is an edge of the form $auv$ and an edge of the form $bwx$ where
$u,v,w,x \in S$. Let $t=|\{u,v\} \cap \{w, x\}|$. If $t=2$, then $a,b,u,v$ forms a copy of $K_4^3$, if $t=3$, then we get a copy of $H_5$ and if $t=4$, then $a,b,u,v,w,x$ forms a copy of $H_6$.      
\qed

Using Theorem \ref{K4H5H6}, Corollary \ref{supercorf}, and the fact that $H_6$ is a subhypergraph of  $K_4^3(2)$ and a subhypergraph of  $H_5(2)$, we obtain the following.

\begin{corollary} \label{cor}
There is an absolute positive constant $c$ such that 
$$r(H_6, K_n^3) < 2^{cn \log n}.$$
\end{corollary}

\subsection{Embedding cycles}
In order to finish the proof, we will show that for all $s \ge 5$, 

(1) $C_8^3 \subset H_6(2)$, 

(2) $C_{s+3}^3 \subset  C_s^3(2)$, and 

(3) $C_{2s}^3 \subset C_s^3(2)$. 

We will apply (1)--(3) together with  Corollaries \ref{supercor} and \ref{cor}. Using (1) we obtain
$$r(C_8^3, K_n^3) \le  r(H_6(2), K_n^3)) <(r(H_6, K_n^3))^{O(1)} < 2^{O(n \log n)}.$$
Using (2) repeatedly we obtain for all $s\equiv 2$ (mod 3) and $s \ge 11$, that
$$r(C_s, K_n^3) \le r(C_8(s), K_n^3) < (r(C_8, K_n^3))^{O(1)} < 2^{O(n \log n)}.$$
Using (3) for $s=16 \equiv 1$ (mod 3) we obtain
$$r(C_{16}, K_n^3) \le  r(C_8(2), K_n^3) <(r(C_8, K_n^3))^{O(1)} < 2^{O(n \log n)}.$$
 Finally, applying (2)  again will finish the proof for all $s\equiv 1$ (mod 3) and $s \ge 16$.
 
 We now turn to the proofs of (1)--(3). Recall that $H_6=\{abu, abv, abw, abx, auv, bwx\}$.
 Showing that $C_8^3 \subset H_6(2)$ is equivalent to producing a homomorphism $\phi$ from $C_8^3$ to $H_6$ where every vertex of $H_6$ has at most two preimages. Define $\phi$ as follows:
 \begin{align*}
 \phi(1)&=u \\
 \phi(2)&=v \\
 \phi(8)=\phi(3)&=a\\
 \phi(7)=\phi(4)&=b\\
 \phi(5)&=x\\
 \phi(6)&=w.
 \end{align*}
 It is easy to check that for every $i \in Z_8$ we have $\phi(i)\phi(i+1)\phi(i+2) \in H_6$.
 Next we show that (2) holds by producing $\phi: Z_{s+3} \rightarrow Z_s$ as follows: let  $\phi(i)=i$ for
 $i \in \{1, \ldots, s\}$, $\phi(s+1)=s-2$, $\phi(s+2)=s-1$ and  $\phi(s+3)=s$.  It is easy to check that $\phi$ is a homomorphism from $V(C_{s+3}^3)$ to $V(C_s^3)$.  Finally, (3) holds due to  $\phi: Z_{2s} \rightarrow Z_s$ as follows: let  $\phi(i)=\phi(s+i)=i$ for
 $i \in \{1, \ldots, s\}$.

\end{document}